\documentstyle[12pt, twoside, leqno]{article}
\input amssym.def
\input amssym.tex

\renewcommand{\theequation}{\arabic{equation}}

\newenvironment{demo}[1]%
{\vskip-\lastskip\medskip
  \noindent
  {\em #1.}\enspace
  }%
{\qed\par\medskip
  }

\newcommand{\qed}{
  \strut\hfill
  \mbox{$\Box$}
  }

\newcommand{\An}{H^{[n]} }
\newcommand{\ch}{\mbox{ch} }
\newcommand{\C}{ \Bbb C }
\newcommand{\D}{ \cal D }
\newcommand{\DA}{{\D}(A)}
\newcommand{\End}{{\rm End}}
\newcommand{\Supp}{{\rm Supp}}
\newcommand{\Fock}{{\cal F}}
\newcommand{\g}{{\gamma}}
\newcommand{\Gh}{{\frak G}_{\hbar}}
\newcommand{\Hx}{{H}}
\newcommand{\hD}{\widehat{\D}}
\newcommand{\hDA}{\widehat{\D}(A)}
\newcommand{\hDH}{\widehat{\D}(\Hx)}
\newcommand{\hf}{\frac12}
\newcommand{\Ln}{L^{[n]}}
\newcommand{\On}{{\cal O}^{[n]}}
\newcommand{\Hn}{H^*(\Xn)}
\newcommand{\Tr}{ {\rm Tr} }
\newcommand{\vac}{|0\rangle}
\newcommand{\W}{{\cal W}_{1+\infty} }
\newcommand{\Xn}{ X^{[n]}}
\newcommand{\Z}{ \Bbb Z }

\newtheorem{theorem}{Theorem}[section]
\newtheorem{lemma}{Lemma}[section]
\newtheorem{remark}{Remark}[section]

\newtheorem{definition}{Definition}[section]

\newtheorem{corollary}{Corollary}[section]

\begin{document}
\title{Stability of the cohomology rings of Hilbert schemes of
points on surfaces}
\author{Wei-Ping Li\thanks{Partially supported by the grant HKUST6170/99P.},
Zhenbo Qin{}\thanks{Partially supported by an NSF grant and an
Alfred P. Sloan Research Fellowship.},
and Weiqiang Wang\thanks{Partially supported by an NSF grant.} }

\date{}
\maketitle

\begin{abstract}
{We establish some remarkable properties of the
cohomology rings of the Hilbert scheme $\Xn$
of $n$ points on a projective surface $X$, from which one sees
to what extent these cohomology rings are (in)dependent of $X$ and $n$.}
\end{abstract}


%\tableofcontents
%%
%%
%%
%%
%%
%%
%%
\section{Introduction}

Lehn \cite{Leh} and more recently the authors \cite{LQW1} have
developed (vertex) algebraic calculus to study the cup products in
the Hilbert schemes $\Xn$ of $n$ points on a projective surface $X$.
This approach
was built on the earlier beautiful formula of G\"ottsche \cite{Got}
on the Betti numbers of $\Xn$ and an important construction of
Heisenberg algebra of Nakajima \cite{Na1, Na2} and
Grojnowski \cite{Gro}.
In \cite{LQW1}, we obtained a set of ring generators for
the rational cohomology ring
$H^*(\Xn)$, which has not been accessible in general by classical
algebro-geometric methods (see however \cite{Mar}).
Using this set of generators, an
algorithm, first pointed out by Lehn \cite{Leh} in a restricted case,
can be given to compute the cup product of any two cohomology classes
in $H^*(\Xn)$ for an {\it arbitrary} projective surface $X$.
This algebraic approach has been surprisingly effective
in establishing new purely geometric results, as indicated
in the further developments of Lehn-Sorger and
the authors \cite{LS1, LQW2, LS2} on the cohomology rings $H^*(\Xn)$.
We refer to \cite{Wa} for a detailed overview and
further references on closely related topics.

In this paper, we establish some remarkable properties of the
cohomology rings of the Hilbert scheme $\Xn$
of $n$-points on a projective surface $X$, from which one sees
to what extent these cohomology rings are (in)dependent of $X$ and $n$.
As a consequence, we are able to introduce
a ring $\frak H_X$ which encodes all the cohomology ring structures of
$H^*(\Xn)$ for all $n$, and further determine its structure.
To achieve these, we will
extensively use and sharpen the techniques developed in
the earlier works of Lehn and of the authors.
%\cite{Leh,LQW1, LQW2} (also cf. \cite{LS2}).
Needless to say, the Heisenberg operators of Nakajima and
Grojnowski are part of our basic vocabulary used in this paper.

We first obtain a quantitative description of
the cup product of the ring generators given in \cite{LQW1},
which indicate to what extent the
cup product of cohomology classes of $\Xn$ is (in)dependent of the
canonical class $K_X$ and Euler class $e_X$ of $X$
(see Theorem~\ref{product_of_G}).
As a corollary (which has been implicit in
the earlier work \cite{Leh, LQW1}),
we see clearly that if there exists a ring isomorphism
 from $H^*(X)$ to $H^*(Y)$ for two projective surfaces $X$ and $Y$
which sends the canonical class $K_X$ to $K_Y$, then the
cohomology rings of the Hilbert schemes $\Xn$ and $Y^{[n]}$ are
isomorphic for any $n$.
In addition, we obtain the general structure of intersection numbers
on $\Xn$ in terms of intersection numbers on $X$.
This general structure bears some similarities with the general
structure of the Donaldson invariants from Donaldson theory
(compare with \cite{EGL}).

Using Theorem~\ref{product_of_G}, we work out the cup products of two
cohomology classes which are monomials of Heisenberg generators,
and observe that the cup products
are independent of $n$ in an appropriate sense
(see Theorem~\ref{cup_product}).
Roughly speaking, Theorem~\ref{cup_product} says that the cup product
of certain cohomology classes in $H^*(\Xn)$ with $n$ being large
can be read off from the cup product of cohomology classes in
$H^*(X^{[m]})$ with $m$ being small.
In other words, the cup product on $X^{[m]}$
partially determines the cup product on $\Xn$ when $n > m$.
This stability result enables us to construct a
super-commutative associative ring ${\frak H}_X$,
called the {\em Hilbert ring} associated to $X$,
which captures {\em all} the information about the cohomology ring
of the Hilbert scheme $\Xn$ for each $n$.
We further prove that $\frak H_X$ is isomorphic to
a super-symmetric algebra with a simple set of generators
which essentially comes from the set of ring generators for
the cohomology rings $H^*(\Xn)$ found in \cite{LQW2}.

In a sequel, we shall develop a counterpart of our results
in terms of the orbifold cohomology rings \cite{CR} of
the symmetric products, and clarify the connections
with our present work. In another direction, it is natural
to expect that results similar to those in the present paper hold
as well for the quantum cohomology rings
of the Hilbert schemes of points on projective surfaces.

The layout of the paper is as follows. In Sect.~\ref{sec_general},
we collect some known results and definitions. In
Sect.~\ref{sec_pushforward}, we establish a series of technical
lemmas related to Heisenberg generators and pushforwards. In
Sect.~\ref{sec_chern}, we work out the cup product of
certain Chern characters in
the cohomology ring of the Hilbert scheme,
and derive some consequences. In Sect.~\ref{sec_stability}, we
establish the stability of the cohomology ring of $\Xn$.
In Sect.~\ref{sec_hilbert}, we introduce the Hilbert ring and
determine its structure.

\medskip\noindent
{\bf Conventions:} All cohomology groups are in $\Bbb Q$-coefficients.
For a continuous map $p: Y_1 \to Y_2$ between two smooth compact
manifolds and for $\alpha_1 \in H^*(Y_1)$, we define $p_*(\alpha_1)$
to be ${\rm PD}^{-1}p_{*}({\rm PD}(\alpha_1))$ where
${\rm PD}$ stands for the Poincar\'e duality.

\bigskip\noindent
{\bf Acknowledgments:} The authors thank the referee for
valuable comments, and Jun Li and Gang Tian for
stimulating discussions. After communicating the paper to
Lehn and Sorger, we were informed that they have also
independently anticipated some of the results in the paper.

\section{Generalities} \label{sec_general}

%In this section, we fix some standard notations and review the
%constructions of several linear operators related to the Hilbert
%schemes of points on surfaces. In addition, we recall some basic
%results from \cite{Na1, Na2, Gro, Leh, LQW1, LQW2}.

Let $X$ be a smooth projective surface over $\C$, and $\Xn$ be the
Hilbert scheme of $n$-points in $X$. An element in the Hilbert scheme
$\Xn$ is represented by a length-$n$ $0$-dimensional closed
subscheme of $X$. It is well-known that $\Xn$ is smooth. Let
${\cal Z}_n=\{(\xi, x) \subset \Xn \times X \, | \, x\in \Supp{(\xi)}\}$,
and $X^n$ be the $n$-th Cartesian product.

\begin{definition} \rm \label{H}
 \begin{enumerate}
 \item[{\rm (i)}] Let $\Bbb H = \oplus_{n, i \ge 0} \Bbb H^{n,i}$ denote
the double graded vector space with $\Bbb H^{n,i} \stackrel{\rm
def}{=} H^i(\Xn)$, and $\Bbb H_n \stackrel{\rm def}{=} H^*(\Xn)
\stackrel{\rm def}{=} \oplus_{i=0}^{4n} H^i(\Xn)$. The element $1$
in $H^0(X^{[0]}) = \Bbb Q$ is called the {\it vacuum vector} and
denoted by $|0\rangle$;

 \item[{\rm (ii)}] $\frak f \in \End(\Bbb H)$ is
{\it homogeneous of bidegree} $(\ell, m)$ if
$\frak f(\Bbb H^{n,i}) \subset \Bbb H^{n+\ell,i+m}$;

 \item[{\rm (iii)}] For $\frak f$ and
$\frak g \in \End(\Bbb H)$ of bidegrees $(\ell, m)$ and $(\ell_1,
m_1)$ respectively, define the {\it Lie superalgebra bracket}
$[\frak f, \frak g]$ by putting $[\frak f, \frak g]
= \frak f \frak g - (-1)^{m m_1} \frak g \frak f$.
 \end{enumerate}
\end{definition}

A non-degenerate super-symmetric bilinear form $(, )$ on $\Bbb H$
is induced from the standard one on $\Bbb H_n=H^*(\Xn)$ defined by
$\displaystyle{(\alpha,\beta) =\int_{\Xn} \alpha\beta}$ for
$\alpha, \beta\in H^*(\Xn)$.
For $\frak f \in \End(\Bbb H)$ of bidegree $(\ell, m)$, we can
define its {\it adjoint} $\frak f^\dagger \in \End(\Bbb H)$ by
$(\frak f(\alpha), \beta) = (-1)^{m \cdot |\alpha|} \cdot
(\alpha, \frak f^\dagger(\beta))$
where $|\alpha| = s$ if $\alpha \in H^s(\Xn)$. Note that the
bidegree of $\frak f^\dagger$ is $(-\ell, m - 4 \ell)$. Also,
for $\frak g \in \End(\Bbb H)$ of bidegree $(\ell_1, m_1)$, we have
\begin{eqnarray} \label{adjoint1}
(\frak f \frak g)^\dagger = (-1)^{m m_1} \cdot \frak g^\dagger
\frak f^\dagger \qquad {\rm and} \qquad [\frak f, \frak g]^\dagger
= -[\frak f^\dagger, \frak g^\dagger].
\end{eqnarray}

We recall that the Heisenberg operators
$\frak a_n(\alpha) \in \End(\Bbb H)$
with $n \in \Z$ and $\alpha \in H^*(X)$ were defined in
\cite{Na1, Gro, Na2}. These operators satisfy the property
$\frak a_n(\alpha) = (-1)^n \cdot \frak a_{-n}(\alpha)^\dagger$.
In the next two sets of definitions, we collect various operators
from \cite{Leh, LQW1, LQW2}.
We have adopted here the usual convention
in the theory of vertex algebras on the signs of indices.
For example, our indices for the Heisenberg and Virasoro generators
coincide with those used in the paper \cite{LS2}, but differ exactly by
a sign from the notations adopted in \cite{Leh, LQW1}.

\begin{definition} \rm \label{operators}
\begin{enumerate}
\item[{\rm (i)}] The {\it normally ordered product}
$:\frak a_{m_1}\frak a_{m_2}:$ is defined by
\begin{eqnarray*}
   :\frak a_{m_1}\frak a_{m_2}:
   & = &
   \left\{
      \everymath{\displaystyle}
      \begin{array}{ll}
  \frak a_{m_1}\frak a_{m_2}, &m_1 \le m_2  \\
  \frak a_{m_2}\frak a_{m_1}, &m_1 \ge m_2.
\end{array}
    \right.
\end{eqnarray*}
For $n \in \Z$, define $\frak L_n: H^*(X) \to \End(\Bbb H)$ by
$\frak L_n = -{1 \over 2} \cdot \sum\limits_{m \in \Z}
:\frak a_m \frak a_{n-m}: \tau_{2*}$.
Here for $k \ge 1$, $\tau_{k*}: H^*(X) \to H^*(X^k)$ is the
map induced by the diagonal embedding $\tau_k: X \to X^k$, and
$\frak a_{m_1} \cdots \frak a_{m_k}(\tau_{k*}(\alpha)) =
\sum_j \frak a_{m_1}(\alpha_{j,1}) \cdots \frak
a_{m_k}(\alpha_{j,k})$
when $\tau_{k*}\alpha = \sum_j \alpha_{j,1} \otimes \cdots \otimes
\alpha_{j, k}$ via the K\"unneth decomposition of $H^*(X^k)$;

\item[{\rm (ii)}] Define the linear operator $\frak d \in \End(\Bbb H)$
by $\frak d = \oplus_{n} c_1(p_{1*}{\cal O}_{{\cal Z}_n})$,
where $p_1$ is the projection of $\Xn \times X$ to $\Xn$,
%$\partial \Xn$ is the boundary of $\Xn$ consisting of
%all $\xi \in \Xn$ with $|\Supp(\xi)| < n$,
and the first Chern class $c_1(p_{1*}{\cal
O}_{{\cal Z}_n})$ of $p_{1*}{\cal O}_{{\cal Z}_n}$ acts on
$\Bbb H_n = H^*(\Xn)$ by the cup product.

\item[{\rm (iii)}] For a linear operator $\frak f \in \End(\Bbb H)$,
define its {\it derivative} $\frak f'$ by
$\frak f' = [\frak d, \frak f]$.
The higher derivative $\frak f^{(k)}$ of $\frak f$ is defined
inductively by $\frak f^{(k)} = [\frak d, \frak f^{(k-1)}]$.
\end{enumerate}
\end{definition}

\begin{definition} \rm \label{GB}
 \begin{enumerate}
 \item[{\rm (i)}] Fix $i, n \ge 0$ and $\alpha \in H^*(X)$.
Let $G_i(\alpha, n)$ denote
the $H^{|\alpha|+2i}(\Xn)$-component of $p_{1*}({\rm ch}({\cal
O}_{{\cal Z}_n}) \cdot p_2^*{\rm td}(X) \cdot p_2^*\alpha) \in
\Bbb H_n$, where $p_1$ and $p_2$ are the two projections of
$\Xn \times X$. Let $B_i(\alpha, n) = 0$ when $i \ge n$,
and $B_i(\alpha, n) = {1/(n-i-1)!} \cdot \frak a_{-1}(1_X)^{n-i-1}
\frak a_{-(i+1)}(\alpha)|0\rangle$ when $i < n$.

 \item[{\rm (ii)}] For $i \ge 0$ and $\alpha \in H^*(X)$,
the {\it Chern character operator}
$\frak G_i(\alpha) \in \End({\Bbb H})$
is defined to be the operator which acts on the component
$\Bbb H_n$ by the cup product by $G_i(\alpha, n)$.
The operator $\frak B_i(\alpha) \in \End({\Bbb H})$
is defined to be the operator which acts on the component
$\Bbb H_n$ by the cup product by $B_i(\alpha, n)$.
\end{enumerate}
\end{definition}

\begin{theorem} \label{commutator}
Let $K_X$ be the canonical divisor of the smooth projective
surface $X$. Let $k \ge 0, n,m \in \Z$ and $\alpha, \beta \in
H^*(X)$. Then,
\begin{enumerate}
\item[{\rm (i)}] $\displaystyle{[\frak a_n(\alpha), \frak a_m(\beta)]
= -n \cdot \delta_{n+m} \cdot \int_X(\alpha \beta) \cdot {\rm Id}_{\Bbb H}}$
where ${\rm Id}_{\Bbb H}$ stands for the identity map of $\Bbb H$,
and $\delta_{n+m}$ is $1$ when $n+m=0$ and $0$ when $n+m \ne 0$;

\item[{\rm (ii)}] $[\frak L_n(\alpha), \frak a_m(\beta)] = -m \cdot \frak
a_{n+m}(\alpha \beta)$;

\item[{\rm (iii)}] $\displaystyle{\frak a_n'(\alpha)
= n \cdot \frak L_n(\alpha) -
{n(|n|-1)/2} \cdot \frak a_n(K_X \alpha)}$;

\item[{\rm (iv)}] $\displaystyle{[\frak G_k(\alpha),
\frak a_{-1}(\beta)] = {1/k!}
\cdot \frak a_{-1}^{(k)}(\alpha \beta)}$;

\item[{\rm (v)}] $\displaystyle{[\ldots [\frak G_k(\alpha),
\frak a_{n_1}(\alpha_1)],
\ldots], \frak a_{n_{k+1}}(\alpha_{k+1})] =
-\prod_{\ell = 1}^{k+1} n_\ell \cdot \frak a_{n_1 + \ldots
+n_{k+1}}(\alpha \alpha_1 \cdots \alpha_{k+1})}$

for all $n_1, \ldots , n_{k+1} \in \Z$ with $\sum_{\ell =1}^{k+1}
n_\ell \ne 0$ and all $\alpha_1, \ldots, \alpha_{k+1} \in H^*(X)$.
\end{enumerate}
\end{theorem}

Theorem \ref{commutator} (i) was proved in \cite{Na2}. The next
two formulas in Theorem \ref{commutator} were obtained in
\cite{Leh}. Theorem \ref{commutator} (iv) and (v) were from
\cite{LQW1}. Also, as observed in \cite{Na1, Gro}, $\Bbb H$
is an irreducible representation of the Heisenberg algebra
generated by the $\frak a_i(\alpha)$'s with the vacuum vector
$|0\rangle \in H^0(X^{[0]})$ being the highest weight vector.
Our next Theorem was proved in \cite{LQW1, LQW2}.

\begin{theorem} \label{basis}
For $n \ge 1$, the cohomology ring $\Bbb H_n = H^*(\Xn)$
is generated by the cohomology classes $G_i(\alpha, n)$
(respectively, the cohomology classes $B_i(\alpha, n)$)
where $0 \le i < n$ and $\alpha$ runs over
a fixed linear basis of $H^*(X)$.
\end{theorem}

\section{Pushforwards and multi-commutators}
\label{sec_pushforward}

In this section, we establish several technical lemmas concerning
the pushforward maps $\tau_{k*}$ and multiple commutators. These
lemmas will be used throughout the paper. We shall also introduce
the concept of {\it a universal linear combination}.

Our first lemma about the pushforward maps $\tau_{k*}$
is elementary but plays an essential role in the
entire paper. We remark that in this lemma and hereafter,
$\tau_{k*}(\alpha)$
is understood to be $\displaystyle{\int_X \alpha}$
when $k = 0$ and $\alpha \in H^*(X)$.

\begin{lemma} \label{pushforward}
Let $k , u \ge 1$ and $\alpha, \beta \in H^*(X)$. Assume that
$\tau_{k*}(\alpha) = \sum_{i} \alpha_{i, 1} \otimes \ldots \otimes
\alpha_{i, k}$ under the K\" unneth decomposition of $H^*(X^k)$.
Then for $0 \le j \le k$, we have
\begin{eqnarray*}
 \tau_{k*}(\alpha\beta)
  &=& \sum_{i} (-1)^{|\beta| \cdot
     \sum_{\ell=j+1}^{k}|\alpha_{i, \ell}|} \cdot
     \left (\otimes_{s=1}^{j-1} \alpha_{i, s} \right )
     \otimes (\alpha_{i, j}\beta)
     \otimes \left ( \otimes_{t=j+1}^{k} \alpha_{i, t} \right )
      \\
 \tau_{(k-1)*}(\alpha \beta)
    &=& \sum_i (-1)^{|\beta| \sum_{\ell=j+1}^k |\alpha_{i,\ell}|}
     \int_X \alpha_{i,j}\beta \cdot
     \otimes_{1 \le s \le k, s \ne j} \alpha_{i, s}
   \\
 \tau_{(k+u-1)*}(\alpha)
  &=&  \sum_{i} \left (\otimes_{s=1}^{j-1} \alpha_{i, s} \right )
     \otimes (\tau_{u*}\alpha_{i, j}) \otimes
     \left ( \otimes_{t=j+1}^{k} \alpha_{i, t} \right ).
\end{eqnarray*}
\end{lemma}
\begin{demo}{Proof}
The basic idea is to use the projection formula. We have
\begin{eqnarray*}
 & &\sum_{i} (-1)^{|\beta| \cdot
     \sum_{\ell=j+1}^{k}|\alpha_{i, \ell}|} \cdot
     \left (\otimes_{s=1}^{j-1} \alpha_{i, s} \right )
     \otimes (\alpha_{i, j}\beta)
     \otimes \left ( \otimes_{t=j+1}^{k} \alpha_{i, t} \right )  \\
 &=&\left (\sum_{i}^{} \alpha_{i, 1} \otimes \ldots
     \otimes \alpha_{i, k} \right ) \cdot p_j^*(\beta)
     = \tau_{k*}(\alpha) \cdot p_j^*(\beta)                 \\
 &=&\tau_{k*}\left (\alpha \cdot (p_j \circ \tau_k)^*(\beta)\right )
     = \tau_{k*}(\alpha  \beta)
\end{eqnarray*}
where $p_j$ is the projection of $X^k$ to the $j$th factor. This
proves the first formula.
The proofs of the second formula and the third formula are similar.
\end{demo}

\begin{lemma} \label{tau_k_tau_{k-1}}
Let $k, s \ge 1$, $n_1, \ldots, n_k, m_1, \ldots, m_s \in \Z$, and
$\alpha, \beta \in H^*(X)$. Then,
 \begin{enumerate}
\item[{\rm (i)}] $[\frak a_{n_1} \cdots \frak a_{n_{k}} (\tau_{k*}\alpha),
\frak a_{m_1} \cdots \frak a_{m_{s}}(\tau_{s*}\beta)]$ is equal to

$$-\sum_{t=1}^k \sum_{j=1}^s n_t \delta_{n_t+m_j} \cdot
\left ( \prod_{\ell=1}^{j-1} \frak a_{m_\ell}
\prod_{1 \le u \le k, u \ne t} \frak a_{n_u} \prod_{\ell=j+1}^{s} \frak
a_{m_\ell} \right )(\tau_{(k+s-2)*}(\alpha\beta));$$

\item[{\rm (ii)}] the derivative $(\frak a_{n_1} \cdots \frak
a_{n_k}(\tau_{k*}\alpha))'$ is equal to
\begin{eqnarray*}
 & &-\sum_{j=1}^k {n_j \over 2} \cdot \sum_{m_1 + m_2 = n_j}
    \frak a_{n_1} \cdots \frak a_{n_{j-1}} :\frak a_{m_1}\frak a_{m_2}:
    \frak a_{n_{j+1}} \cdots \frak a_{n_k}(\tau_{(k+1)*}\alpha)   \\
 & &- \sum_{j=1}^k {n_j(|n_j|-1) \over 2} \cdot
    \frak a_{n_1} \cdots \frak a_{n_k}(\tau_{k*}(K_X\alpha)).
\end{eqnarray*}
 \end{enumerate}
\end{lemma}
\begin{demo}{Proof}
Follows from Theorem~\ref{commutator} (i) and (iii),
and Lemma~\ref{pushforward}.
\end{demo}

\begin{definition} \rm \label{universal}
Let $X$ be a projective surface, $s \ge 1$, and $\alpha_1,
\ldots, \alpha_s \in H^*(X)$. Let $k_1, \ldots, k_s \ge 0$, and
$n_{i,j} \in \Z$ with $1 \le i \le s$ and $1 \le j \le k_i$. Then,
a {\it universal linear combination} of $\frak a_{n_{i, 1}} \cdots
\frak a_{n_{i, k_i}}(\tau_{k_i*}(\alpha_i))$, $1 \le i \le s$ is a
linear combination of the form $\displaystyle{\sum_{i=1}^s
f_i(k_i,n_{i, 1}, \ldots, n_{i, k_i}) \frak a_{n_{i, 1}} \cdots
\frak a_{n_{i, k_i}}(\tau_{k_i*}(\alpha_i))}$ where the
coefficients $f_i(k_i,n_{i, 1}, \ldots, n_{i, k_i})$ are
independent of $X, \alpha_1, \ldots, \alpha_s$. A {\it universal
linear combination} of $\frak a_{n_{i, 1}} \cdots \frak a_{n_{i,
k_i}} (\tau_{k_i*}(\alpha_i))|0\rangle$, $1 \le i \le s$ is
defined in a similar way.
\end{definition}

\begin{lemma} \label{derivative}
Let $k, s \ge 0$, $n, m_1, \ldots, m_s \in \Z$, and $\alpha,
\beta_1, \ldots, \beta_s \in H^*(X)$. Then,
 \begin{enumerate}
\item[{\rm (i)}] $\frak a_n^{(k)}(\alpha)$ is a universal
linear combination of
$\frak a_{n_1} \cdots \frak a_{n_{k-r+1}}
(\tau_{(k-r+1)*}(K_X^r \alpha))$
where $0 \le r \le 2$ and $n_1+\ldots+n_{k-r+1} = n$;

\item[{\rm (ii)}] $[\cdots [\frak a_n^{(k)}(\alpha), \frak
a_{m_{1}}(\beta_1)], \cdots], \frak a_{m_{s}}(\beta_s)]$ is a
universal linear combination of
\begin{eqnarray*}
\frak a_{n_1} \cdots \frak a_{n_{k-s-r+1}}
 (\tau_{(k-s-r+1)*}(K_X^r \alpha\beta_1 \cdots \beta_{s}))
\end{eqnarray*}
where $0 \le r \le 2$ and $n_1+\ldots+n_{k-s-r+1} = n + m_1+\ldots +m_{s}$.
 \end{enumerate}
\end{lemma}
\begin{demo}{Proof}
Since $K_X^3 = 0$, (i) follows from repeatedly applying Lemma
\ref{tau_k_tau_{k-1}} (ii). Now (ii) follows from (i) and
repeatedly applying Lemma \ref{tau_k_tau_{k-1}} (i).
\end{demo}

\begin{lemma} \label{switch}
Let $e_X$ denote the Euler class of $X$. Fix $k \ge 2$,
$n_1, \ldots, n_k \in \Z$, and $\alpha \in H^*(X)$.
Let $j$ satisfy $1 \le j < k$. Then, $\frak a_{n_1} \cdots
\frak a_{n_k}(\tau_{k*}\alpha)$ is equal to
\begin{eqnarray*}
\left ( \prod_{1 \le s < j} \frak a_{n_s} \cdot \frak
a_{n_{j+1}} \frak a_{n_{j}} \cdot \prod_{j+1 < s \le k} \frak
a_{n_s} \right ) (\tau_{k*}\alpha) - n_j \delta_{n_j+n_{j+1}}
\prod_{1 \le s \le k \atop s \ne j, j+1} \frak
a_{n_s}(\tau_{(k-2)*}(e_X\alpha)).
\end{eqnarray*}
\end{lemma}
\begin{demo}{Proof}
Note that $\sum_t \beta_{t, 1}\beta_{t, 2} = e_X \beta$ if
$\tau_{2*}(\beta) = \sum_t \beta_{t, 1} \otimes \beta_{t, 2}$.
Now our result follows from Theorem~\ref{commutator}~(i)
and the second and third formulas in Lemma~\ref{pushforward}.
\end{demo}

\begin{lemma} \label{f_{12}}
Let $n \ge 1$, $\alpha \in H^*(X)$, and $\frak f \in
\End({\Bbb H})$ with $\frak f'=0$. Then,
\begin{eqnarray*}
[\frak f, \frak a_{-(n+1)}(\alpha)] =-{1 \over n} \cdot \{ [[\frak
f, \frak a_{-1}(1_X)]', \frak a_{-n}(\alpha)] + [\frak
a_{-1}'(1_X), [\frak f, \frak a_{-n}(\alpha)]]\}.
\end{eqnarray*}
\end{lemma}
\begin{demo}{Proof}
Appeared implicitly in \cite{Leh},
and follows from Theorem~\ref{commutator}~(ii), (iii).
\end{demo}

\section{Products of Chern characters} \label{sec_chern}

In this section, we prove that the products of Chern characters
$G_k(\alpha, n)$ can be written as some universal finite linear
combination of monomials of Heisenberg generators
(see Theorem~\ref{product_of_G} below). As an application,
we obtain the general structure of intersection numbers on
the Hilbert scheme $\Xn$. We remark that Theorem~\ref{product_of_G}
will also be used substantially in later sections.

The following lemma is a variation of the Lemma~5.26 in \cite{LQW1}.

\begin{lemma} \label{nonsense1}
Fix $k \ge 0$ and $b \ge 1$. Let $\frak g \in \End(\Bbb H)$ be of
bidegree $(\tilde s, s)$ satisfying
\begin{eqnarray} \label{nonsense1.1}
[[\cdots [\frak g, \frak a_{n_1}(\alpha_1)], \cdots], \frak
a_{n_{k+1}}(\alpha_{k+1})] =0
\end{eqnarray}
for any $n_1, \ldots, n_{k+1} < 0$ and $\alpha_{1}, \ldots,
\alpha_{k+1} \in H^*(X)$. Let $A = \frak a_{m_1}(\beta_{1}) \cdots
\frak a_{m_b}(\beta_{b})|0\rangle$ where $m_1, \ldots, m_b < 0$
and $\beta_{1}, \ldots, \beta_b \in H^*(X)$. Then, $\frak g(A)$ is
equal to
\begin{eqnarray*}
&&\sum_{i=0}^{k} \sum_{\sigma_i}
  (-1)^{s \sum\limits_{\ell \in \sigma_i^0} |\beta_\ell|
  + \sum\limits_{j=1}^i \sum\limits_{\ell \in \sigma_i^0, \ell > \sigma_i(j)}
  |\beta_{\sigma_i(j)}||\beta_\ell|} \cdot  \\
&&\cdot \prod_{\ell \in \sigma_i^0} \frak a_{m_\ell}(\beta_{\ell})
 [[\cdots [\frak g, \frak a_{m_{\sigma_i(1)}}(\beta_{\sigma_i(1)})],
 \cdots], \frak a_{m_{\sigma_i(i)}}(\beta_{\sigma_i(i)})] |0\rangle
\end{eqnarray*}
where for each fixed $i$, $\sigma_i$ runs over all the maps
$\{\, 1, \ldots, i \,\} \to \{\, 1, \ldots, b\,\}$
satisfying $\sigma_i(1) < \cdots < \sigma_i(i)$, and
$\sigma_i^0 = \{ \ell \, | \, 1 \le \ell \le b, \ell \ne
\sigma_i(1), \ldots, \sigma_i(i) \}$.    \qed
\end{lemma}

\begin{lemma} \label{G_k}
Let $s \ge 1$, and $\alpha, \beta \in H^*(X)$. Then,
$[\frak G_k(\alpha), \frak a_{n_1} \cdots \frak a_{n_s}(\tau_{s*}\beta)]$
is a universal linear combination of expressions
$\frak a_{m_1} \cdots \frak a_{m_{k+s-r}}
(\tau_{(k+s-r)*}(K_X^r\alpha\beta))$
where $0 \le r \le 2$ and $m_1+\ldots+m_{k+s-r} = n_1+\ldots+n_{s}$.
\end{lemma}
\begin{demo}{Proof}
First of all, let $s = 1$. Note that $\frak a_{0}(\beta) = 0$.
Since $\frak G_k(\alpha)^\dagger = \frak G_k(\alpha)$ and
$\frak a_{n_1}(\beta)^\dagger = (-1)^{n_1} \frak a_{-n_1}(\beta)$,
we see from (\ref{adjoint1}) that
$[\frak G_k(\alpha), \frak a_{n_1}(\beta)]^\dagger =
-[\frak G_k(\alpha)^\dagger, \frak a_{n_1}(\beta)^\dagger]$
$= (-1)^{1+n_1}[\frak G_k(\alpha), \frak a_{-n_1}(\beta)]$.
Since $(\frak a_{m_1} \cdots \frak a_{m_{k+s-r}}
(\tau_{(k+s-r)*}(K_X^r\alpha\beta)))^\dagger$ is equal to
\begin{eqnarray*}
(-1)^{m_1+ \ldots + m_{k+s-r}} \cdot
\frak a_{-m_{k+s-r}} \cdots \frak a_{-m_1}
(\tau_{(k+s-r)*}(K_X^r\alpha\beta)),
\end{eqnarray*}
we need only to prove the statement for $[\frak G_k(\alpha), \frak
a_{n_1}(\beta)]$ with $n_1 \le -1$. When $n_1 = -1$, $[\frak
G_k(\alpha), \frak a_{n_1}(\beta)] =1/k! \cdot \frak
a_{-1}^{(k)}(\alpha\beta)$. So the statement for $s=1$ and $n_1 =
-1$ follows from Lemma \ref{derivative} (i). When $n_1 \le -2$, we
see from Lemma \ref{f_{12}} that
\begin{eqnarray*}
 & &[\frak G_k(\alpha), \frak a_{n_1}(\beta)]   \\
 &=&{1 \over n_1+1} \cdot \{ [[\frak G_k(\alpha), \frak a_{-1}(1_X)]',
    \frak a_{n_1+1}(\beta)] + [\frak a_{-1}'(1_X),
    [\frak G_k(\alpha), \frak a_{n_1+1}(\beta)]]\}  \\
 &=&{1 \over n_1+1} \cdot \left \{ \right .
    [[\frak G_k(\alpha), \frak a_{-1}(1_X)]',
    \frak a_{n_1+1}(\beta)]    \\
 & &\quad + [\frak a_{-1}(1_X),
    [\frak G_k(\alpha), \frak a_{n_1+1}(\beta)]]' -
    [\frak a_{-1}(1_X), [\frak G_k(\alpha), \frak a_{n_1+1}(\beta)]']
    \left . \right \}.
\end{eqnarray*}
So the statement for $s = 1$ and $n_1 \le -2$ follows from
induction and Lemma \ref{tau_k_tau_{k-1}}.

Next, let $s \ge 2$. Let $\tau_{s*}(\beta) = \sum_{i}
\beta_{i, 1} \otimes \ldots \otimes \beta_{i, s} \in H^*(X^s)$.
Then, we have
$[\frak G_k(\alpha), \frak a_{n_1} \cdots \frak a_{n_s}
    (\tau_{s*}\beta)]
=\sum_{i} [\frak G_k(\alpha), \frak a_{n_1}(\beta_{i, 1}) \cdots
    \frak a_{n_s}(\beta_{i, s})]$.
By symmetry, it suffices to show that
$\sum_{i} [\frak G_k(\alpha), \frak a_{n_1}(\beta_{i, 1})] \frak
a_{n_2}(\beta_{i, 2}) \cdots \frak a_{n_s}(\beta_{i, s})$
is a universal linear combination of the forms
$\frak a_{m_1} \cdots \frak
a_{m_{k+s-r}}(\tau_{(k+s-r)*}(K_X^r\alpha\beta))$
where $0 \le r \le 2$ and $m_1+\ldots+m_{k+s-r} =
n_1+\ldots+n_{s}$. To prove this, we apply what we have already
proved in the preceding paragraph to $[\frak G_k(\alpha), \frak
a_{n_1}(\beta_{i, 1})]$. So $[\frak G_k(\alpha),
\frak a_{n_1}(\beta_{i, 1})]$ equals
\begin{eqnarray*}
\sum_{0 \le r \le 2 \atop m_1+\ldots+m_{k-r+1} = n_1} f_r(k, n_1,
m_1, \ldots, m_{k-r+1}) \frak a_{m_1} \cdots \frak a_{m_{k-r+1}}
(\tau_{(k-r+1)*}(K_X^r \alpha\beta_{i, 1}))
\end{eqnarray*}
where $f_r(k, n_1, m_1, \ldots, m_{k-r+1})$ stands for universal
rational numbers independent of $X$, $\alpha$, and $\beta_{i, 1}$.
In particular, these universal numbers are independent of $i$. So
\begin{eqnarray*}
 & &\sum_{i} [\frak G_k(\alpha), \frak a_{n_1}(\beta_{i, 1})]
    \frak a_{n_2}(\beta_{i, 2}) \cdots \frak a_{n_s}(\beta_{i, s}) \\
 &=&\sum_{0 \le r \le 2 \atop m_1+\ldots+m_{k-r+1} = n_1}
    f_r(k, n_1, m_1, \ldots, m_{k-r+1}) \cdot \\
 & &\cdot \sum_{i} \left ( \frak a_{m_1} \cdots \frak a_{m_{k-r+1}}
    (\tau_{(k-r+1)*}(K_X^r\alpha\beta_{i, 1})) \right )
    \frak a_{n_2}(\beta_{i, 2}) \cdots \frak a_{n_s}(\beta_{i, s}).
\end{eqnarray*}
By the first formula in Lemma \ref{pushforward},
$\sum_i (K_X^r\alpha\beta_{i, 1}) \otimes \beta_{i, 2} \otimes
\cdots \otimes \beta_{i, s} = \tau_{s*}(K_X^r\alpha\beta).$
So by the third formula in Lemma \ref{pushforward},
$\sum_i \tau_{(k-r+1)*}(K_X^r\alpha\beta_{i, 1}) \otimes \beta_{i, 2}
\otimes \cdots \otimes \beta_{i, s} = \tau_{(k+s-r)*}(K_X^r\alpha\beta)$.
It follows that $\sum_{i} [\frak G_k(\alpha), \frak
a_{n_1}(\beta_{i, 1})] \frak a_{n_2}(\beta_{i, 2}) \cdots \frak
a_{n_s}(\beta_{i, s})$ equals
\begin{eqnarray*}
\sum_{0 \le r \le 2 \atop m_1+\ldots+m_{k-r+1} = n_1} f_r(k, n_1,
m_1, \ldots, m_{k-r+1})  \left ( \prod_{\ell=1}^{k-r+1} \frak
a_{m_\ell}  \prod_{t=2}^{s} \frak a_{n_t} \right
)(\tau_{(k+s-r)*}(K_X^r\alpha\beta))
\end{eqnarray*}
where $f_r(k, n_1, m_1, \ldots, m_{k-r+1})$ are independent of
$X$, $\alpha$, and $\beta$.
\end{demo}

\begin{definition}   \rm   \label{a_{k,j}}
\begin{enumerate}
\item[{\rm (i)}] Let $s \ge 1$, and $\alpha_1, \ldots, \alpha_s
\in H^*(X)$ be homogeneous.
For a partition $\pi =\{ \pi_1, \ldots, \pi_j \}$
of the set $\{1, \ldots, s \}$, we fix the orders of the elements
listed in each subset $\pi_i$ ($1 \le i \le j$)
once and for all, and define $\ell(\pi) = j$,
$\alpha_{\pi_i} ={\prod_{m \in \pi_i} \alpha_m}$,
and ${\rm sign}(\alpha, \pi)$ by
$\prod_{i=1}^j \alpha_{\pi_i} = {\rm sign}(\alpha, \pi) \cdot
\prod_{i=1}^s \alpha_i.$

\item[{\rm (ii)}] We denote ${\bf 1}_{-n} = 1/n! \cdot
\frak a_{-1}(1_X)^n$ when $n \ge 0$, and ${\bf 1}_{-n}
= 0$ when $n < 0$.
\end{enumerate}
\end{definition}

The geometric meaning of ${\bf 1}_{-n}$ is that
${\bf 1}_{-n}|0\rangle = 1_{X^{[n]}}$ (the fundamental
class of $X^{[n]}$). Also, the choice of the orders for the elements
listed in each $\pi_i$, $1 \le i \le \ell(\pi)$
will affect ${\rm sign}(\alpha, \pi)$,
but will not affect the expression (\ref{product_of_G1})
in our next Theorem where an empty product
$\displaystyle{\prod_{i} G_{k_i}(\alpha_i, n)}$ stands for
$1_{X^{[n]}} = {\bf 1}_{-n}|0\rangle$ by convention.

\begin{theorem} \label{product_of_G}
Let $n \ge 1, s \ge 0$, $k_1, \ldots, k_s \ge 0$, and $\alpha_1,
\ldots, \alpha_s \in H^*(X)$ be homogeneous.
Then, $\displaystyle{\prod_{i=1}^s G_{k_i}(\alpha_i, n)}$
is a finite linear combination of expressions:
\begin{eqnarray} \label{product_of_G1}
\qquad {\rm sign}(\alpha, \pi) \cdot {\bf 1}_{-\left
(n-\sum\limits_{i=1}^{\ell(\pi)}
  \sum\limits_{j=1}^{m_i-r_i} n_{i, j} \right )} \left (
\prod_{i = 1}^{\ell(\pi)} \left ( \prod_{j = 1}^{m_i-r_i} \frak
a_{-n_{i, j}} \right ) (\tau_{(m_i-r_i)*}(\epsilon_i
\alpha_{\pi_i})) \right ) |0\rangle
\end{eqnarray}
whose coefficients  are
independent of  $X$,
$\alpha_1, \ldots, \alpha_s$, and the integer $n$.
Here $\pi$ runs over all partitions of $\{1, \ldots, s
\}$, $\epsilon_i \in \{1_X, K_X, K_X^2, e_X\}$,
$r_i = |\epsilon_i|/2 \le m_i \le 2+ \sum_{j \in \pi_i} k_j,$
$0 < n_{i, 1} \le \ldots \le n_{i, m_i-r_i}$,
$\sum\limits_{j=1}^{m_i-r_i}
n_{i, j} \le \sum\limits_{j \in \pi_i} (k_j+1)$ for each $i$, and
\begin{eqnarray} \label{product_of_G1.1}
\sum_{i = 1}^{\ell(\pi)} \left ( m_i - 2 + \sum_{j = 1}^{m_i-r_i}
n_{i, j} \right ) = \sum_{i=1}^s k_i.
\end{eqnarray}
\end{theorem}
\begin{demo}{Proof}
Use induction on $s$. When $s=0$, the statement is trivial by
our convention. Next, let $s \ge 1$. By induction,
$\displaystyle{\prod_{i=2}^s G_{k_i}(\alpha_i, n)}$
is a linear combination of expressions:
\begin{eqnarray*}
& &{\rm sign}(\alpha, \sigma) \cdot
   {\bf 1}_{-(n-\tilde n)}
   \left ( \prod_{i = 1}^{\ell(\sigma)}
   \left ( \prod_{j = 1}^{m_i-r_i} \frak a_{-n_{i, j}} \right )
   (\tau_{(m_i-r_i)*}(\epsilon_i \alpha_{\sigma_i}))
   \right ) |0\rangle    \\
&=&{{\rm sign}(\alpha, \sigma) \over (n-\tilde n)!} \cdot
   \frak a_{-1}(1_X)^{n-\tilde n} \left ( \prod_{i = 1}^{\ell(\sigma)}
   \left ( \prod_{j = 1}^{m_i-r_i} \frak a_{-n_{i, j}} \right )
   (\tau_{(m_i-r_i)*}(\epsilon_i \alpha_{\sigma_i}))
   \right ) |0\rangle
\end{eqnarray*}
where $\sigma$ runs over all partitions of $\{2, \ldots, s
\}$, $\epsilon_i \in \{1_X, K_X, K_X^2, e_X\}$,
$r_i = |\epsilon_i|/2 \le m_i \le 2+\sum_{j \in \sigma_i} k_j,$
$0 < n_{i, 1} \le \ldots \le n_{i, m_i-r_i}$,
$\sum\limits_{j=1}^{m_i-r_i}  n_{i, j} \le \sum\limits_{j \in
\sigma_i} (k_j+1)$, and $\tilde n = \sum\limits_{i=1}^{\ell(\sigma)}
\sum\limits_{j=1}^{m_i-r_i} n_{i, j}$. Moreover, the coefficients
in the linear combination are independent of $X, \alpha_2, \ldots,
\alpha_s$ and $n$. Now apply $\frak G_{k_1}(\alpha_1)$
to $\displaystyle{\prod_{i=2}^s G_{k_i}(\alpha_i, n)}$,
and move $\frak G_{k_1}(\alpha_1)$ to the
right by using Lemma \ref{nonsense1}.
Note that $|\tau_{(m_i-r_i)*}(\epsilon_i
\alpha_{\sigma_i})| \equiv |\alpha_{\sigma_i}| \pmod 2$.
By Theorem~\ref{commutator}~(v),
$[\ldots [\frak G_{k_1}(\alpha_1),
\frak a_{\ell_1}(\beta_1)], \ldots], \frak
a_{\ell_{k_1+2}}(\beta_{k_1+2})]=0$ when
$\ell_1, \ldots, \ell_{k_1+2} < 0$.
Since $\frak G_{k_1}(\alpha_1)|0\rangle = 0$ and
$\displaystyle{\prod_{i=1}^s G_{k_i}(\alpha_i, n)
= \frak G_{k_1}(\alpha_1) \left (
\prod_{i=2}^s G_{k_i}(\alpha_i, n) \right ),}$
we see from Lemma \ref{nonsense1} that
$\displaystyle{\prod_{i=1}^s G_{k_i}(\alpha_i, n)}$ is
a universal linear combination of expressions:
\begin{eqnarray}   \label{product_of_G2}
 &&\quad {{\rm sign}(\alpha, \sigma) \over (n-\tilde n)!}
   {n-\tilde n \choose t}
   \cdot (-1)^{|\alpha_1|\sum\limits_{1 \le i \le \ell(\sigma), i \not\in U}
   |\alpha_{\sigma_i}|+
   \sum\limits_{v=1}^u \sum\limits_{w > i_v, w \not \in U}
   |\alpha_{\sigma_{i_v}}||\alpha_{\sigma_w}|}  \cdot  \\
 &&\cdot \frak a_{-1}(1_X)^{n-\tilde n -t}
   \left ( \prod_{1 \le i \le \ell(\sigma) \atop i \not\in U}
   \left ( \prod_{j = 1}^{m_i-r_i} \frak a_{-n_{i, j}} \right )
   (\tau_{(m_i-r_i)*}(\epsilon_i \alpha_{\sigma_i}))
   \right )  \cdot  \nonumber \\
 &&\cdot \left [\cdots [\frak G_{k_1}(\alpha_1),
    \underbrace{\frak a_{-1}(1_X)],
    \cdots], \frak a_{-1}(1_X)]}_{t ~\rm{times}},
    \left ( \prod_{j = 1}^{m_{i_1}-r_{i_1}} \frak a_{-n_{{i_1}, j}} \right )
    (\tau_{(m_{i_1}-r_{i_1})*}(\epsilon_{i_1} \alpha_{\sigma_{i_1}}))],
    \right . \nonumber    \\
 &&\left . \quad\quad\quad \cdots ],
    \left ( \prod_{j = 1}^{m_{i_u}-r_{i_u}} \frak a_{-n_{{i_u}, j}} \right )
    (\tau_{(m_{i_u}-r_{i_u})*}(\epsilon_{i_u} \alpha_{\sigma_{i_u}}))
    \right ] |0\rangle    \nonumber
\end{eqnarray}
where $0 \le t \le (k_1+1), u \ge 0$, $(t+u) \ge 1$, $U = \{i_1,
\ldots, i_u\} \subset \{1, \ldots, \ell(\sigma) \}$ with
$i_1 < \ldots < i_u$. Let $\pi$ be the partition of $\{1,
\ldots, s\}$ consisting of all the $\sigma_i$ with $1 \le i \le
\ell(\sigma)$ and $i \not\in U$, and $\displaystyle{\{ 1 \} \coprod
\left ( \coprod_{i \in U} \sigma_i \right )}$. Then,
$\alpha_{\pi_{\ell(\pi)}}
=\alpha_1 \alpha_{\sigma_{i_1}} \cdots \alpha_{\sigma_{i_u}}$ and
\begin{eqnarray} \label{product_of_G3}
{\rm sign}(\alpha, \pi) ={\rm sign}(\alpha, \sigma) \cdot
(-1)^{|\alpha_1|\sum\limits_{1 \le i \le \ell(\sigma), i \not\in U}
|\alpha_{\sigma_i}|+ \sum\limits_{v=1}^u \sum\limits_{w > i_v, w
\not \in U} |\alpha_{\sigma_{i_v}}||\alpha_{\sigma_w}|}.
\end{eqnarray}
In view of (\ref{product_of_G2}) and (\ref{product_of_G3}),
$\displaystyle{\prod_{i=1}^s G_{k_i}(\alpha_i, n)}$ is a linear
combination of
\begin{eqnarray}   \label{product_of_G4}
 &&\quad {\rm sign}(\alpha, \pi)
   \cdot {\bf 1}_{-(n - \tilde n -t)}
   \left ( \prod_{1 \le i < \ell(\pi)}
   \left ( \prod_{j = 1}^{m_i-r_i} \frak a_{-n_{i, j}} \right )
   (\tau_{(m_i-r_i)*}(\epsilon_i \alpha_{\pi_i}))
   \right )  \cdot  \\
 &&\cdot \left [\cdots [\frak G_{k_1}(\alpha_1),
    \underbrace{\frak a_{-1}(1_X)],
    \cdots], \frak a_{-1}(1_X)]}_{t ~\rm{times}},
    \left ( \prod_{j = 1}^{m_{i_1}-r_{i_1}} \frak a_{-n_{{i_1}, j}} \right )
    (\tau_{(m_{i_1}-r_{i_1})*}(\epsilon_{i_1} \alpha_{\sigma_{i_1}}))],
    \right . \nonumber    \\
 &&\left . \quad\quad\quad \cdots ],
    \left ( \prod_{j = 1}^{m_{i_u}-r_{i_u}} \frak a_{-n_{{i_u}, j}} \right )
    (\tau_{(m_{i_u}-r_{i_u})*}(\epsilon_{i_u} \alpha_{\sigma_{i_u}}))
    \right ] |0\rangle.    \nonumber
\end{eqnarray}
with all the coefficients being independent of $X, \alpha_1,
\ldots, \alpha_s$ and $n$. Also notice that in the expression
(\ref{product_of_G4}), the only factor depending on $n$ is ${\bf
1}_{-(n - \tilde n -t)}$.

Let $t = 0$. Then, $u \ge 1$. By Lemma~\ref{G_k},
Lemma \ref{tau_k_tau_{k-1}} (i) and Lemma~\ref{switch},
each expression (\ref{product_of_G4}) is a
universal linear combination of expressions of the form
\begin{eqnarray*}
& &{\rm sign}(\alpha, \pi)
   \cdot {\bf 1}_{-(n - \tilde n)}
   \left ( \prod_{1 \le i < \ell(\pi)}
   \left ( \prod_{j = 1}^{m_i-r_i} \frak a_{-n_{i, j}} \right )
   (\tau_{(m_i-r_i)*}(\epsilon_i \alpha_{\pi_i}))
   \right )  \cdot  \\
& &\quad \cdot \frak a_{-n_{1}} \cdots \frak a_{-n_{m-r}}
   (\tau_{(m-r)*}(\epsilon \alpha_1 \alpha_{\sigma_{i_1}}
    \cdots \alpha_{\sigma_{i_u}}))|0\rangle   \\
&=&{\rm sign}(\alpha, \pi)
   \cdot {\bf 1}_{-(n - \tilde n)}
   \left ( \prod_{1 \le i < \ell(\pi)}
   \left ( \prod_{j = 1}^{m_i-r_i} \frak a_{-n_{i, j}} \right )
   (\tau_{(m_i-r_i)*}(\epsilon_i \alpha_{\pi_i}))
   \right )  \cdot  \\
& &\quad \cdot \frak a_{-n_{1}} \cdots \frak a_{-n_{m-r}}
   (\tau_{(m-r)*}(\epsilon \alpha_{\pi_{\ell(\pi)}}))|0\rangle
\end{eqnarray*}
which is of the form (\ref{product_of_G1}).
Here $\epsilon = {\tilde \epsilon}\epsilon_{i_1} \cdots \epsilon_{i_u}$
with ${\tilde \epsilon} \in \{1_X, K_X, K_X^2, e_X\}$,
$r = |{\tilde \epsilon}|/2 + \sum_{j=1}^u r_{i_j},$
$m = k_1 + \sum_{j=1}^u m_{i_j} - 2(u-1)$,
$0 < n_1 \le \ldots \le n_{m-r}$, and
\begin{eqnarray*}
& &n_1 + \ldots + n_{m-r} = \sum_{i \in U}
   \sum_{j=1}^{m_i-r_i} n_{i, j} \le
   \sum_{i \in U} \sum_{j \in \sigma_i} (k_j+1)   \\
&<&(k_1+1) + \sum_{i \in U} \sum_{j \in \sigma_i} (k_j+1)
   = \sum_{j \in \pi_{\ell(\pi)}} (k_j+1).
\end{eqnarray*}
Note that either $\epsilon = 0$ or $\epsilon \in \{1_X, K_X, K_X^2, e_X\}$.
When $\epsilon \in \{1_X, K_X, K_X^2, e_X\}$,
we have $r = |\epsilon|/2 \le m$.
Since $m_{i_v} \le 2+\sum_{j \in \sigma_{i_v}} k_j$ for
$1 \le v \le u$, we obtain
$$m = k_1 + m_{i_1} + \ldots + m_{i_u} -2(u-1)
\le 2+\sum_{j \in \pi_{\ell(\pi)}} k_j.$$

Next, assume that $t \ge 1$. Then by Theorem~\ref{commutator}
(iv), we have
$$[\cdots [\frak G_{k_1}(\alpha_1),
\underbrace{\frak a_{-1}(1_X)], \cdots], \frak a_{-1}(1_X)]}_{t
~\rm{times}} = {1 \over k_1!} \cdot [\cdots [\frak
a_{-1}^{(k_1)}(\alpha_1), \underbrace{\frak a_{-1}(1_X)], \cdots],
\frak a_{-1}(1_X)]}_{(t-1) ~\rm{times}}$$ which by Lemma
\ref{derivative} (ii), is a universal linear combination of
expressions of the form
$\frak a_{n_1} \cdots \frak a_{n_{k_1-t+2-{\tilde r}}}
(\tau_{(k_1-t+2-{\tilde r})*}(K_X^{\tilde r} \alpha_1))$
where $0 \le {\tilde r} \le 2$ and
$n_1+\ldots+n_{k_1-t+2-{\tilde r}} = -t$. So by
Lemma~\ref{tau_k_tau_{k-1}}~(i) and Lemma~\ref{switch},
(\ref{product_of_G4}) is a universal linear combination of
\begin{eqnarray*}
& &{\rm sign}(\alpha, \pi)
   \cdot {\bf 1}_{-(n - \tilde n -t)}
   \left ( \prod_{1 \le i < \ell(\pi)}
   \left ( \prod_{j = 1}^{m_i-r_i} \frak a_{-n_{i, j}} \right )
   (\tau_{(m_i-r_i)*}(\epsilon_i \alpha_{\pi_i}))
   \right )  \cdot  \\
& &\quad \cdot \frak a_{-n_{1}} \cdots \frak a_{-n_{m-r}}
   (\tau_{(m-r)*}(\epsilon \alpha_1 \alpha_{\sigma_{i_1}}
    \cdots \alpha_{\sigma_{i_u}}))|0\rangle   \\
&=&{\rm sign}(\alpha, \pi)
   \cdot {\bf 1}_{-(n - \tilde n -t)}
   \left ( \prod_{1 \le i < \ell(\pi)}
   \left ( \prod_{j = 1}^{m_i-r_i} \frak a_{-n_{i, j}} \right )
   (\tau_{(m_i-r_i)*}(\epsilon_i \alpha_{\pi_i}))
   \right )  \cdot  \\
& &\quad \cdot \frak a_{-n_{1}} \cdots \frak a_{-n_{m-r}}
   (\tau_{(m-r)*}(\epsilon \alpha_{\pi_{\ell(\pi)}}))|0\rangle
\end{eqnarray*}
which again is of the form (\ref{product_of_G1}). Here $m = (k_1-t+2)
+ m_{i_1} + \ldots + m_{i_u}-2u < 2+ \sum_{j \in
\pi_{\ell(\pi)}} k_j$, $\epsilon \in \{1_X, K_X, K_X^2, e_X\}$, $r
=|\epsilon|/2 \le m$, $0 < n_1 \le \ldots \le n_{m-r}$, and $n_1 +
\ldots + n_{m-r} = t+\sum_{i \in U} \sum_{j=1}^{m_i-r_i} n_{i, j}
\le \sum_{j \in \pi_{\ell(\pi)}} (k_j+1)$ since $t \le (k_1+1)$.

Finally, the cohomology degree of
$\displaystyle{\prod_{i=1}^s G_{k_i}(\alpha_i, n)}$ is equal to
$\sum_{i=1}^s (2k_i+|\alpha_i|)$. Comparing this with the
cohomology degree of (\ref{product_of_G1}), we obtain
(\ref{product_of_G1.1}).
\end{demo}

\begin{corollary} \label{structure}
Let $X$ and $Y$ be two complex projective surfaces. Assume that
there exists a ring isomorphism $\Phi: H^*(X) \to H^*(Y)$ with
$\Phi(K_X) = K_Y$. Then for every $n \ge 1$, the two cohomology
rings $H^*(\Xn)$ and $H^*(Y^{[n]})$ are isomorphic.
\end{corollary}
\begin{demo}{Proof}
Note that $\Phi(e_X) = e_Y$. Since the Chern characters
$G_k(\alpha, n)$ generate the cohomology ring $H^*(\Xn)$,
our result follows from Theorem~\ref{product_of_G}.
\end{demo}

Next, we apply Theorem~\ref{product_of_G} to study intersection
numbers in the Hilbert scheme $\Xn$. For this purpose,
we establish the notation
$\displaystyle{\langle w \rangle = \int_Y w}$
where $w \in H^*(Y)$ and $Y$ stands for a smooth projective variety.

\begin{corollary} \label{intersection_number}
Let $n, s \ge 1$, $k_1, \ldots, k_s \ge 0$, and let $\alpha_1,
\ldots, \alpha_s \in H^*(X)$ be homogeneous cohomology classes.
Assume $\sum\limits_{i=1}^s (2k_i + |\alpha_i|) = 4n$.
Then, $\displaystyle{\left \langle \prod_{i=1}^s
G_{k_i}(\alpha_i, n) \right \rangle}$ is a finite linear combination of
$\displaystyle{\,\, {\rm sign}(\alpha, \pi) \cdot \prod_{i = 1}^{\ell(\pi)}
\langle \epsilon_i \alpha_{\pi_i} \rangle}$
where $\pi$ runs over all partitions of $\{1, \ldots, s
\}$, $\epsilon_i \in \{1_X, K_X, K_X^2, e_X\}$.
Moreover, all the coefficients in this linear combination are
independent of $X$, $\alpha_1, \ldots, \alpha_s$ and $n$.
\end{corollary}
\begin{demo}{Proof}
Note that the positive generator of $H^{4n}(\Xn) \cong \Bbb Q$ is
$\frak a_{-1}([x])^n |0\rangle$ where $[x] \in H^4(X)$ stands for
the cohomology class corresponding to a point $x \in X$.
So by Theorem~\ref{product_of_G},
an expression (\ref{product_of_G1}) nontrivially contributing to
$\displaystyle{\left \langle \prod_{i=1}^s
G_{k_i}(\alpha_i, n) \right \rangle}$ must satisfy:
i) $n_{i, j} = 1$ for all $1 \le i \le \ell(\pi)$ and $1 \le j \le m_i - r_i$;
ii) $\epsilon_i \alpha_{\pi_i} = \langle \epsilon_i \alpha_{\pi_i}
\rangle \cdot [x]$ for all $1 \le i \le \ell(\pi)$;
and iii) $n-\sum_{i=1}^{\ell(\pi)} (m_i-r_i) =0$.
Since $\tau_{k*}([x]) =
\underbrace{[x] \otimes \cdots \otimes [x]}_{k \, {\rm times}}$
for all $k \ge 0$, our conclusion
follows immediately from (\ref{product_of_G1}).
\end{demo}

\section{The stability} \label{sec_stability}

In this section, we establish a remarkable stability for the
cohomology rings of the Hilbert schemes of $n$-points on
projective surfaces as $n$ varies.

We need two lemmas which sharpen the Lemma~3.20
and Lemma~3.5 in \cite{LQW2}. In the first lemma,
we determine the {\em leading} monomial of
Heisenberg generators in
$\displaystyle{\prod_{i=1}^s G_{k_i}(\alpha_i, n)}$.
In the second lemma, we express
$\displaystyle{{\bf 1}_{-\left (n - \sum_{i=1}^s n_i \right )}
\left ( \prod\limits_{i=1}^{s} \frak a_{-n_{i}}(\alpha_{i})
\right ) |0\rangle}$ as a universal finite linear combination of
cup products of the form $\displaystyle{\prod\limits_{j=1}^t
G_{m_j}(\beta_j, n)}$.

\begin{lemma} \label{leading-term}
Let notations be the same as in Theorem~\ref{product_of_G}.
\begin{enumerate}
\item[{\rm (i)}] An expression of the form (\ref{product_of_G1})
satisfying $\sum\limits_{i=1}^{\ell(\pi)} \sum\limits_{j=1}^{m_i-r_i}
n_{i, j} = \sum\limits_{i=1}^s (k_i+1)$
is equal to $\displaystyle{{\bf 1}_{-(n - n_0)}
\prod\limits_{i=1}^{s} \frak a_{-(k_{i}+1)}(\alpha_{i}) \cdot
|0\rangle}$, where $n_0 \stackrel{\rm def}{=} \sum_{i=1}^s (k_i+1)$.

\item[{\rm (ii)}] The coefficient of $\displaystyle{{\bf 1}_{-(n - n_0)}
\prod\limits_{i=1}^{s} \frak a_{-(k_{i}+1)}(\alpha_{i}) \cdot \vac}$
in $\displaystyle{\prod_{i=1}^s G_{k_i}(\alpha_i, n)}$ is
$\displaystyle{\prod_{i=1}^s {(-1)^{k_i} \over (k_i+1)!}}.$
 \end{enumerate}
\end{lemma}
\begin{demo}{Proof}
(i) We may let $s \ge 1$. Since $\displaystyle{\sum\limits_{i=1}^{\ell(\pi)}
\sum\limits_{j=1}^{m_i-r_i} n_{i, j} = \sum_{i=1}^s (k_i+1)}$,
we see from Theorem~\ref{product_of_G}
that $\sum\limits_{j=1}^{m_i-r_i} n_{i, j} = \sum\limits_{j \in \pi_i}
(k_j+1)$ for every $i$. By (\ref{product_of_G1.1}), we obtain
\begin{eqnarray} \label{leading_term1}
\sum_{i = 1}^{\ell(\pi)} (m_i - 2 + |\pi_i|) = 0
\end{eqnarray}
where $|\pi_i|$ stands for the number of elements in the subset $\pi_i$.
Note that for every $i$ with $1 \le i \le \ell(\pi)$,
we have $(m_i-r_i) \ge 1$ since
$\displaystyle{\sum\limits_{j=1}^{m_i-r_i} n_{i, j} =
\sum\limits_{j \in \pi_i} (k_j+1) \ge 1.}$
So for every $i$, $m_i \ge 1$, and $r_i = 0$ if $m_i = 1$.
By (\ref{leading_term1}), $m_i = |\pi_i|=1$ for $1 \le i \le \ell(\pi)$.
Thus for $1 \le i \le \ell(\pi)$, we have $r_i = 0$, $\epsilon_i = 1_X$,
and $n_{i, 1} = (k_j+1)$ if $\pi_i = \{ j \}$.

Now let $\pi_i = \{ t_i \}$ for $1 \le i \le \ell(\pi) = s$. Then the
expression (\ref{product_of_G1}) is
\begin{eqnarray*}
{\rm sign}(\alpha, \pi) \cdot {\bf 1}_{-(n-n_0)} \left (
   \prod_{i = 1}^{s} \frak a_{-(k_{t_i}+1)}(\alpha_{t_i})
   \right ) |0\rangle
={\bf 1}_{-(n - n_0)} \left ( \prod\limits_{i=1}^{s}
   \frak a_{-(k_{i}+1)}(\alpha_{i}) \right ) |0\rangle.
\end{eqnarray*}

(ii) The idea is to use induction on $s$ and track the proof of
Theorem~\ref{product_of_G} more carefully.
When $s = 0$, the statement is trivial.
Next, let $s \ge 1$ and ${\tilde n}_0 = \sum\limits_{j=2}^s
(k_j+1)$. Assume that the coefficient of
${\bf 1}_{-(n - {\tilde n}_0)} \left ( \prod\limits_{i=2}^{s}
\frak a_{-n_{i}}(\alpha_{i}) \right ) |0\rangle$ in the cup product
$\displaystyle{\prod_{i=2}^s G_{k_i}(\alpha_i, n)}$ is equal to
$\displaystyle{\prod_{i=2}^s {(-1)^{k_i} \over (k_i+1)!}}$.
Tracking the proof of Theorem~\ref{product_of_G} and applying
Theorem \ref{commutator} (v), we conclude that the coefficient of
$\displaystyle{{\bf 1}_{-(n - n_0)} \left ( \prod\limits_{i=1}^{s}
\frak a_{-(k_{i}+1)}(\alpha_{i})\right ) |0\rangle}$ in
$\displaystyle{\prod_{i=1}^s G_{k_i}(\alpha_i, n)}$ is equal to
$\displaystyle{\prod_{i=1}^s {(-1)^{k_i} \over (k_i+1)!}}$.
\end{demo}

\begin{lemma} \label{combination_of_G}
Fix $n, s \ge 1$, $n_1, \ldots, n_s \ge 1$, and $\alpha_1,
\ldots, \alpha_s\in H^*(X)$. Put $n_0 = \sum\limits_{i=1}^s n_i$.
Then, $\displaystyle{{\bf 1}_{-(n - n_0)} \left (
\prod\limits_{i=1}^{s} \frak a_{-n_{i}}(\alpha_{i}) \right )
|0\rangle}$ is a finite linear combination of
\begin{eqnarray} \label{combination_of_G1}
\prod\limits_{j=1}^t G_{m_j}(\beta_j, n)
\end{eqnarray}
whose coefficients are independent of $X,
\alpha_1, \ldots, \alpha_s$ and $n$.
Here $\sum\limits_{j=1}^t (m_j+1) \le n_0$, and $\beta_1, \ldots,
\beta_t$ depend only on $e_X, K_X, \alpha_1, \ldots, \alpha_s$ and
$\tau_{i*}$ with $1 \le i \le n_0$.
\end{lemma}
\begin{demo}{Proof}
We use induction on $n_0$. When $n_0 = 1$, $s=n_1 = 1$. By
the Lemma~3.20~(i) in \cite{LQW2}, ${\bf 1}_{-(n - 1)} \frak
a_{-1}(\alpha_{1})|0\rangle = G_0(\alpha_1, n)$. So the lemma
holds for $n_0 = 1$.

Next, let $n_0 > 1$. Let $k_i = n_i -1$. Then,
$k_i \ge 0$ for every $i$. By Theorem~\ref{product_of_G},
$\displaystyle{\prod_{i=1}^s G_{k_i}(\alpha_i, n)}$
is a finite linear combination of expressions of the form
(\ref{product_of_G1}) such that the coefficients in this linear
combination are independent of $X, \alpha_1, \ldots, \alpha_s$ and $n$.
Note that $\sum\limits_{i=1}^{\ell(\pi)}
\sum\limits_{j=1}^{m_i-r_i} n_{i, j} \le
\sum\limits_{i=1}^{\ell(\pi)} \sum\limits_{j \in \pi_i} (k_j+1)
= \sum_{i=1}^s (k_i +1) = n_0$. By induction, those expressions
(\ref{product_of_G1}) with $\sum\limits_{i=1}^{\ell(\pi)}
\sum\limits_{j=1}^{m_i-r_i} n_{i, j} < n_0$ are linear
combinations of the form (\ref{combination_of_G1}) where
$\sum\limits_{j=1}^t (m_j+1) \le (n_0-1)$, and $\beta_1, \ldots,
\beta_t$ depend only on $e_X, K_X, \alpha_1, \ldots, \alpha_s,
\tau_{i*}$ with $1 \le i \le (n_0-1)$. Moreover, the coefficients
in these linear combinations are independent of $X, \alpha_1,
\ldots, \alpha_s$ and $n$. Now our lemma follows from
Lemma~\ref{leading-term}.
\end{demo}

\begin{remark} \label{leading-term-in-combination_of_G}
{\rm By Lemma \ref{combination_of_G},
an expression (\ref{combination_of_G1}) in
$\displaystyle{{\bf 1}_{-(n - n_0)} \left ( \prod\limits_{i=1}^{s}
\frak a_{-n_{i}}(\alpha_{i}) \right )|0\rangle}$ satisfies
$\sum\limits_{j=1}^t (m_j+1) \le n_0$. In fact, we see from the
proof of Lemma \ref{combination_of_G} that
an expression (\ref{combination_of_G1})
satisfies the upper bound $\sum\limits_{j=1}^t (m_j+1) = n_0$ if
and only if it is equal to $\displaystyle{\prod_{i=1}^s
G_{n_i-1}(\alpha_i, n)}$ whose coefficient  is
$\displaystyle{\prod_{i=1}^s ((-1)^{n_i-1} n_i!)}$
in view of Lemma~\ref{leading-term} (ii).}
\end{remark}

Next, we prove a lemma which says that the Chern character
$G_k(\alpha, n)$ can be expressed as a universal finite linear
combination of cup products
$\displaystyle{\prod\limits_{j=1}^t B_{m_j}(\beta_j, n)}$
(see Definition \ref{GB}~(i)).
In other words, our next lemma essentially reverses the process
in Lemma \ref{combination_of_G}. This lemma will be used later
in the proof of Theorem \ref{th_structure}.

\begin{lemma} \label{combination_of_B}
The Chern character $G_k(\alpha, n)$ is a finite linear
combination of products
$\displaystyle{\prod\limits_{j=1}^t B_{m_j}(\beta_j, n)}$
whose coefficients   are
independent of $X$, $\alpha$ and $n$.
Here $\sum\limits_{j=1}^t m_j \le k$, and
$\beta_1, \ldots, \beta_t$
depend only on $e_X$, $K_X$ and $\alpha$.
In addition,  $\sum\limits_{j=1}^t m_j = k$
if and only if the product $\displaystyle{\prod\limits_{j=1}^t
B_{m_j}(\beta_j, n)}$  equals $B_k(\alpha, n)$ whose coefficient
 is $(-1)^k/(k+1)!$.
\end{lemma}
\begin{demo}{Proof}
Use induction on $k$. When $k=0$, we have $G_0(\alpha, n) =
B_0(\alpha, n)$ by the Lemma 3.20 (i) in \cite{LQW2}.
Next, we assume that the lemma is true for $0,
\ldots, k-1$ for some fixed $k \ge 1$. We shall prove that the
lemma holds for $k$ as well. We apply Lemma~\ref{combination_of_G}
and Remark~\ref{leading-term-in-combination_of_G} to
$B_k(\alpha, n) = {\bf 1}_{-(n-k-1)} \frak a_{-(k+1)}(\alpha)|0\rangle.$
We see that $\displaystyle{G_k(\alpha, n) - {(-1)^k \over (k+1)!} \cdot
B_k(\alpha, n)}$ is a finite linear combination of
$\displaystyle{\prod\limits_{j=1}^u G_{n_j}(\gamma_j, n)}$
where $\sum\limits_{j=1}^u (n_j+1) < (k+1)$, and $\gamma_1,
\ldots, \gamma_u$ depend only on $e_X, K_X$ and $\alpha$.
Moreover, the coefficients in this linear combination are
independent of $X, \alpha$ and $n$. Note that $n_j < k$ for all $1
\le j \le u$. So by induction hypothesis, the lemma holds for $k$.
\end{demo}

\begin{remark} \label{newproof}
{\rm Lemma~\ref{combination_of_G} and Lemma~\ref{combination_of_B}
provide a new proof to Theorem~\ref{basis}
which was originally proved in \cite{LQW1, LQW2}.
}
\end{remark}

Our stability result below indicates that the cup product on
the Hilbert scheme $\Xn$ are independent of $n$
in an appropriate sense. Furthermore, we find an explicit form of
the leading term in the cup product. This result enables us
to construct a ring and determine its structure in the next section.

\begin{theorem} \label{cup_product}
Let $s \ge 1$ and $k_i \ge 1$ for $1 \le i \le s$. Fix $n_{i, j}
\ge 1$ and $\alpha_{i, j} \in H^*(X)$ for $1 \le j \le k_i$, and
fix $n$ with $n \ge \sum\limits_{j=1}^{k_i} n_{i, j}$ for all $1 \le i
\le s$. Then the cup product
\begin{eqnarray} \label{cup_product1}
\prod_{i=1}^s \left ( {\bf 1}_{-(n - \sum_{j=1}^{k_i} n_{i, j})}
\left (\prod_{j=1}^{k_i} \frak a_{-n_{i, j}}(\alpha_{i, j}) \right
) |0\rangle \right )
\end{eqnarray}
in $H^*(\Xn)$ is equal to a finite linear combination of
monomials of the form
\begin{eqnarray} \label{cup_product2}
 {\bf 1}_{-(n - \sum_{p=1}^N m_{p})}
\left ( \prod_{p=1}^N \frak a_{-m_{p}}(\g_{p}) \right ) |0\rangle
\end{eqnarray}
whose coefficients are independent of $X,
\alpha_{i,j}$ and $n$. Here $\sum\limits_{p=1}^N m_{p} \le
\sum\limits_{i=1}^s \sum\limits_{j=1}^{k_i}
n_{i,j}$, and $\g_1, \ldots, \g_N$ depend only on $e_X, K_X,
\alpha_{i,j}$, $1\le i \le s, 1\le j\le k_i$.  In addition, the expression
(\ref{cup_product2}) satisfies the upper bound
$\sum\limits_{p=1}^N m_{p} = \sum\limits_{i=1}^s
\sum\limits_{j=1}^{k_i} n_{i,j}$ if and only if it is equal to
$\displaystyle{{\bf 1}_{-(n -\sum_{i=1}^s \sum_{j=1}^{k_i} n_{i,j})}
\left ( \prod_{i=1}^s\prod_{j=1}^{k_i} \frak a_{-n_{i,
j}}(\alpha_{i, j}) \right ) |0\rangle}$ whose coefficient is $1$.
\end{theorem}
\begin{demo}{Proof}
Put $N_{i} = \sum_{j=1}^{k_i} n_{i, j}$ for $1 \le i \le s$.
For each $i$, we see from Lemma~\ref{combination_of_G} that
$\displaystyle{{\bf 1}_{-(n - \sum_{j=1}^{k_i} n_{i, j})}
\left (\prod_{j=1}^{k_i}
\frak a_{-n_{i, j}}(\alpha_{i, j}) \right ) |0\rangle}$
is a finite linear combination of products
$\displaystyle{\prod\limits_{j=1}^{t_i} G_{m_{i,j}}(\beta_{i,j}, n)}$
where $\sum\limits_{j=1}^{t_i} (m_{i,j}+1) \le N_{i}$, and
$\beta_{i,1}, \ldots, \beta_{i,t_i}$ depend only on $e_X, K_X,$ $
\alpha_{i,1}, \ldots, \alpha_{i,k_i}$ and
$\tau_{j*}$ with $1 \le j \le N_{i}$. Moreover, the
coefficients in the linear combinations are independent of $X,
\alpha_{i,1}, \ldots, \alpha_{i,k_i}$ and $n$. By Remark
\ref{leading-term-in-combination_of_G}, the product satisfies the
upper bound $\sum\limits_{j=1}^{t_i} (m_{i,j}+1) = N_{i}$ if
and only if it is equal to $\displaystyle{\prod_{j=1}^{k_i}
G_{n_{i,j}-1}(\alpha_{i,j}, n)}$. Furthermore, the coefficient of
$\displaystyle{\prod_{j=1}^{k_i} G_{n_{i,j}-1}(\alpha_{i,j}, n)}$
in $\displaystyle{{\bf 1}_{-(n - \sum_{j=1}^{k_i} n_{i, j})} \left
( \prod_{j=1}^{k_i} \frak a_{-n_{i, j}}(\alpha_{i, j}) \right )
|0\rangle}$ is equal to
$\displaystyle{\prod_{j=1}^{k_i}((-1)^{n_{i,j}-1} n_{i,j}!)}$.

So (\ref{cup_product1}) is a universal finite linear combination
of products $\displaystyle{\prod_{i=1}^s \prod\limits_{j=1}^{t_i}
G_{m_{i,j}}(\beta_{i,j}, n)}$
where $\sum\limits_{i=1}^s \sum\limits_{j=1}^{t_i} (m_{i,j}+1) \le
\sum\limits_{i=1}^s N_{i} =\sum\limits_{i=1}^s
\sum\limits_{j=1}^{k_i} n_{i,j}$. Also,
$\sum\limits_{i=1}^s \sum\limits_{j=1}^{t_i} (m_{i,j}+1)
=\sum\limits_{i=1}^s \sum\limits_{j=1}^{k_i} n_{i,j}$
if and only if the product $\displaystyle{\prod_{i=1}^s
\prod\limits_{j=1}^{t_i} G_{m_{i,j}}(\beta_{i,j}, n)}$
is equal to $\displaystyle{\prod_{i=1}^s
\prod_{j=1}^{k_i} G_{n_{i,j}-1}(\alpha_{i,j}, n)}$. The
coefficient of $\displaystyle{\prod_{i=1}^s \prod_{j=1}^{k_i}
G_{n_{i,j}-1}(\alpha_{i,j}, n)}$ in (\ref{cup_product1}) is equal to
$\displaystyle{\prod_{i=1}^s \prod_{j=1}^{k_i}((-1)^{n_{i,j}-1}
n_{i,j}!)}$.

It follows from Theorem~\ref{product_of_G} that (\ref{cup_product1})
is a universal linear combination of expressions
(\ref{cup_product2}). The statement for the expression
(\ref{cup_product2}) reaching the upper bound $\sum_{p=1}^N m_{p}
= \sum_{i=1}^s \sum_{j=1}^{k_i} n_{i,j}$ follows from Lemma
\ref{leading-term}.
\end{demo}

\section{The Hilbert ring} \label{sec_hilbert}

Using the stability result proved in the previous section,
we shall introduce and determine {\it the Hilbert
ring} ${\frak H}_X$ associated to a projective surface $X$.

Given a finite set $S$ which is a disjoint union of subsets $S_0$
and $S_1$, we denote by ${\mathcal P}(S)$ the set of
partition-valued functions $\rho =(\rho(c))_{c \in S}$ on $S$ such
that for every $c \in S_1$, the partition $\rho(c)$ is required to
be {\it strict} in the sense that $\rho(c) =(1^{m_1(c)} 2^{m_2(c)}
\ldots )$ with $m_r(c) = 0$ or $1$ for all $r \ge 1$.

Now let us take a linear basis $S= S_0 \cup S_1$ of $H^*(X)$ such
that $1_X \in S_0$, $S_0 \subset H^{\rm even}(X)$ and $S_1 \subset
H^{\rm odd}(X)$. If we write $\rho =(\rho (c))_{c \in S}$ and
$\rho(c) =(r^{m_r(c)})_{r \ge 1} =(1^{m_1(c)} 2^{m_2(c)} \ldots)$,
then we put $\displaystyle{\ell(\rho) = \sum_{c \in S} \ell(\rho (c))
  = \sum_{c\in S, r\geq 1} m_r(c)}$ and
\begin{eqnarray*}
\Vert \rho \Vert = \sum_{c \in S} |\rho (c)|
  =\sum_{c\in S, r\geq 1} r \cdot m_r(c),  \quad
{\mathcal P}_n(S) = \{\rho \in{\mathcal P} (S)\;|\;\; \Vert \rho
\Vert =n\}.
\end{eqnarray*}
Given $\rho=(\rho(c))_{c\in S}=(r^{m_r(c)})_{c\in S,r \ge 1} \in
{\mathcal P}(S)$ and $n \ge 0$, we define
\begin{eqnarray*}
  {\frak a}_{- \rho(c)}(c)
  &=& \prod_{r \ge 1}
  {\frak a}_{-r}(c)^{m_r(c)} = {\frak a}_{-1}(c)^{m_1(c)}
  {\frak a}_{-2}(c)^{m_2(c)} \cdots    \\
 {\frak a}_{\rho}(n)
 &=& {\bf 1}_{-(n-\Vert \rho \Vert)}
  \prod_{c\in S}{\frak a}_{-{\rho}(c)}(c) \cdot \vac \in H^*(\Xn)
\end{eqnarray*}
where we fix the order of the elements $c \in S_1$ appearing in
$\displaystyle{\prod_{c \in S}}$ once and for all.
Note from Definition \ref{a_{k,j}} (ii) that ${\frak a}_{\rho}(n)=0$
for $0 \le n < \Vert \rho \Vert$.

As $\rho$ runs over all partition-valued functions on $S$ with
$\Vert \rho \Vert \le n$, the corresponding ${\frak a}_{\rho}(n)$
linearly span $H^*(\Xn)$ as a corollary to the theorem of Nakajima
and Grojnowski \cite{Na2}. By Theorem ~\ref{cup_product} (for $s=2$),
we have the cup product
\begin{eqnarray} \label{eq_structure}
{\frak a}_{\rho}(n) \cdot {\frak a}_{\sigma}(n) =
\sum_{\nu}d_{\rho\sigma}^\nu {\frak a}_{\nu}(n)
\end{eqnarray}
in $H^*(\Xn)$, where $\Vert \nu \Vert \le
\Vert \rho \Vert +\Vert \sigma \Vert$
and the structure coefficients $d_{\rho\sigma}^\nu$ are
independent of $n$. Even though the cohomology classes
${\frak a}_{\nu}(n)$ with $\Vert \nu \Vert \le n$ in $H^*(\Xn)$
are not linearly independent in general, we have the following.

\begin{lemma} \label{lem_unique}
The structure constants $d_{\rho\sigma}^\nu$ in the formula
(\ref{eq_structure}) are uniquely determined by the requirement
that they are independent of $n$.
\end{lemma}
\begin{demo}{Proof}
Assume that there exist a finite subset $I \subset {\mathcal P}(S)$
and some constants $c^\nu \in \Bbb Q$ independent of $n$ such that
for all $n \ge 0$, we have
\begin{eqnarray} \label{eq_indep}
\sum_{\nu \in I}c^\nu {\frak a}_{\nu}(n)=0.
\end{eqnarray}
As an immediate consequence of the theorem of
Nakajima and Grojnowski \cite{Na2}, the Heisenberg monomials
$\prod_{c\in S}\frak a_{-\rho(c)}(c) \cdot \vac$,
where $\rho =(\rho(c))_c \in {\cal P}_n(S)$, are linearly
independent in the cohomology ring $H^*(\Xn)$. Therefore,
by the definition of $\frak a_{\nu}(n)$, we may assume
in (\ref{eq_indep}) that any two distinct $\nu$ and $\tilde \nu$
in $I$ satisfy $\tilde \nu(c) =\nu(c)$ for $c \neq 1_X$, $\nu(1_X)
=(1^{m_1}2^{m_2} \cdots)$ and
$\tilde \nu(1_X) =(1^{m_1+\Vert \tilde \nu \vert-\Vert \nu \Vert}
2^{m_2}\cdots)$
(here we assume for definiteness that
$n \ge \Vert \tilde \nu\Vert > \Vert \nu\Vert$). In this case,
we have ${\frak a}_{\nu}(n)= {(n-\Vert \tilde \nu\Vert)!/
(n-\Vert\nu\Vert)!} \cdot {\frak a}_{\tilde \nu}(n)$.
Letting $n \to \infty$, we see from (\ref{eq_indep})
that $c^{\tilde \nu}=0$ for the $\tilde \nu \in I$ with
the largest size $\Vert \tilde \nu\Vert$.
So all the constants $c^\nu$ are zero.
\end{demo}

Now we are ready to introduce the Hilbert ring.

\begin{definition}  \rm
The {\it Hilbert ring} associated to a projective surface $X$,
denoted by ${\frak H}_X$, is defined to be the ring with a linear
basis formed by the symbols ${\frak a}_\rho$, $\rho \in {\mathcal
P}(S)$ and with the multiplication defined by
${\frak a}_{\rho} \cdot {\frak a}_{\sigma} =
\sum_{\nu} d_{\rho\sigma}^\nu {\frak a}_{\nu}$
where the structure constants $d_{\rho\sigma}^\nu$ are from
the relations (\ref{eq_structure}).
\end{definition}

Note that the Hilbert ring does not depend on the choice of a
linear basis $S$ of $H^*(X)$ containing $1_X$ since the operator
${\frak a}_n(\alpha)$ depends on the cohomology class
$\alpha \in H^*(X)$ linearly.
It follows from the super-commutativity and associativity of the
cohomology ring $H^*(\Xn)$ that the Hilbert ring ${\frak H}_X$
itself is also super-commutative and associative.
The ring ${\frak H}_X$ captures all the information of
the cohomology rings of $\Xn$ for all $n$,
as we easily recover the relations (\ref{eq_structure}) from the
ring ${\frak H}_X$. We summarize these observations into the following.

\begin{theorem} {\rm (Stability)} \label{th_stab}
For a given projective surface $X$, the cohomology rings
$H^*(\Xn)$, $n\ge 1$ give rise to a Hilbert ring ${\frak H}_X$
which completely encodes the cohomology ring structure of
$H^*(\Xn)$ for each $n$.
\end{theorem}

We further have the following result on the structure of the
Hilbert ring ${\frak H}_X$. For convenience,
in the case when $\ell(\rho)=1$, that is, when the partition
$\rho(c)$ is a one-part partition $(r)$ for some element $c\in S$
and is empty for all the other elements in $S$, we will simply write
${\frak a}_{\rho} ={\frak a}_{r,c}$ and
${\frak a}_{\rho}(n) ={\frak a}_{r,c}(n)$.

\begin{theorem} \label{th_structure}
The Hilbert ring $\frak H_X$ is isomorphic to the tensor product
$P \otimes E$, where $P$ is the polynomial algebra generated
by $\frak a_{r,c},\;c\in S_0, r \ge 1$ and $E$ is the exterior
algebra generated by ${\frak a}_{r,c},\;c\in S_1, r \ge 1$.
\end{theorem}
\begin{demo}{Proof}
Note that ${\frak a}_{r,c}(n)= {\bf 1}_{-(n-r)}
\frak a_{-r}(c)\vac = B_{r-1}(c, n)$.
By Lemmas~\ref{combination_of_G} and
\ref{combination_of_B}, the ring ${\frak H}_X$ is
generated by the elements ${\frak a}_{r,c},$ where
$c\in S =S_0 \cup S_1$ and $r \ge 1$.

By the super-commutativity of $\frak H_X$, we have ${\frak
a}_{r,c}^2 =0$ for $c\in S_1$ and $r \ge 1$. It remains to show
that as $\rho =(r^{m_r(c)})_{s \in S, r \ge 1}$ runs over
${\mathcal P}(S)$, the monomials
$\displaystyle{\prod_{c\in S, r\ge 1}{\frak a}_{r, c}^{m_r(c)}}$
are linearly independent in ${\frak H}_X$. Assume
$\displaystyle{\sum_{i \in I} d_i \prod_{c\in S, r\ge 1}
{\frak a}_{r, c}^{m^i_r(c)} =0}$
where $d_i \in \Bbb Q$ and $\rho_i=(r^{m^i_r(c)})_{c \in S, r \ge
1}$ runs over a finite set $I$ of distinct elements in ${\mathcal
P}(S)$. By the definition of the structure
constants in ${\frak H}_X$ and ${\bf 1}_{-(n-r)}
\frak a_{-r}(c)\vac = {\frak a}_{r,c}(n)$,
\begin{eqnarray} \label{eq_class}
 \quad \sum_{i \in I} d_i \cdot \prod_{c\in S, r\ge 1}
 \left ( {\bf 1}_{-(n-r)}
 \frak a_{-r}(c)\vac \right )^{m^i_r(c)}
=\sum_{i \in I} d_i \cdot \prod_{c\in S, r\ge 1}
 \frak a_{r, c}(n)^{m^i_r(c)}
=0.
\end{eqnarray}

Take an integer $n$ large enough such that
$n \ge n_i \stackrel{\rm def}{=} \sum_{r,c}r m^i_{r}(c)$
for all $i \in I$. By Theorem~\ref{cup_product},
Eq.~(\ref{eq_class}) can be rewritten as
\begin{eqnarray} \label{eq_stab3}
\sum_{i \in I} d_i \left ( {\bf 1}_{-(n-n_i)} \prod_{c \in S,
r\ge 1} (\frak a_{-r}(c))^{m^i_r(c)} \cdot \vac + w_i \right ) = 0
\end{eqnarray}
where each $w_i$ is a finite linear combination of
${\bf 1}_{-(n - \sum_{p=1}^N m_{p})}
\prod_{p=1}^N \frak a_{-m_{p}}(\g_{p}) \cdot |0\rangle$
%expressions of the form
%\begin{eqnarray*}
%{\bf 1}_{-(n - \sum_{p=1}^N m_{p})}
%\prod_{p=1}^N \frak a_{-m_{p}}(\g_{p}) \cdot |0\rangle
%\end{eqnarray*}
with $\sum_{p=1}^N m_{p} < n_i$ and $\g_p \in S$ for every $p$.
Recall that ${\bf 1}_{-k} =1/{k!} \cdot {\frak a}_{-1}(1_X)^k$, $k
\ge 0$. If we multiply (\ref{eq_stab3}) by $n!$, and locate in the
resulting summation those terms whose coefficients contain
the largest power of $n$, then we see that
\begin{eqnarray} \label{eq_class1}
\sum_i d_i \cdot {\bf 1}_{-(n - n_i)} \prod_{c \in S, r\ge 1}
(\frak a_{-r}(c))^{m^i_r(c)} \cdot \vac =0
\end{eqnarray}
where $i$ satisfies $n_i = {\rm max}\{ n_j| j \in I \}$.
Since all the integers $n_i$ in (\ref{eq_class1}) are equal,
the Heisenberg monomials in (\ref{eq_class1}) are linearly independent
as a corollary to the theorem of Nakajima and Grojnowski \cite{Na2}.
Thus all the coefficients $d_i$ in (\ref{eq_class1}) are zero.
By repeating the above argument,
we obtain that $d_i = 0$ for all $i \in I$.
\end{demo}

\noindent Department of Mathematics, HKUST, Clear Water Bay,
Kowloon, Hong Kong, mawpli@uxmail.ust.hk; \\
Department of Mathematics, University of Missouri, Columbia, MO
65211, USA, zq@math.missouri.edu;\\
Department of Mathematics, University of Virginia,
Charlottesville, VA 22904, ww9c@virginia.edu.
\end{document}